\documentclass{article}
\usepackage{graphicx}
\usepackage{amsmath}

\input{tcilatex}

\begin{document}

\title{\textbf{The Generalized Method of Exhaustion}}
\author{Anthony A. Ruffa \\
Naval Undersea Warfare Center Division\\
1176 Howell Street\\
Newport, RI 02841-1708}
\maketitle

\begin{abstract}
The method of exhaustion is generalized to a simple formula that can be used
to integrate functions under very general conditions, provided that the
integral exists. Both a geometric proof (following the usual procedure for
the method of exhaustion) and an independent algebraic proof are provided.
Applications provided as examples include use of the formula to generate new
series for common functions as well as computing the group velocity
distribution resulting from waves diffracted from an aperture.
\end{abstract}

\section{Introduction}

\noindent The method of exhaustion was used by the ancient Greeks to
integrate a limited number of functions before the development of calculus
[1]. It will be shown here that this ancient method can be generalized to a
simple formula to integrate functions under very general conditions,
assuming that the integral exists over the interval in question. The formula
is

\bigskip

\begin{equation}
\int_{a}^{b}f(x)dx=(b-a)\sum_{n=1}^{\infty
}\sum_{m=1}^{2^{n}-1}(-1)^{m+1}2^{-n}f(a+m(b-a)/2^{n})
\end{equation}

\bigskip

The proof is simple and geometric in nature, directly following the
procedure originally used for the method of exhaustion. An independent
algebraic proof is also provided for the special case for which a = 0.

\section{Geometric Proof}

The geometric proof is a generalization of the successive approximation
method historically used to obtain the exact area under the simple curves
previously integrated by this approach. Consider an arbitrary function f(x)
that is piecewise continuous on [a,b]. The procedure to find the area under
f(x) on [a,b] involves successive approximations with triangles that first
intersect f(x) at the midpoint, x = a + (b - a)/2, and then at the new
midpoints, x = a + (b - a)/4, x = a + 3(b - a)/4, and so on. The first such
approximation is shown in figure 1. Its area is

\begin{equation}
A_{1}=\frac{1}{2}(b-a)f(a+(b-a)/2).
\end{equation}

\bigskip 

The second such approximation adds two new triangles (figure 2), each
sharing two vertices with the original triangle, and each bisecting one of
the remaining arcs of f(x) at its third vertex. The two triangles have the
area

\begin{equation}
A_{2}=\frac{1}{4}(b-a)\left\{
f(a+(b-a)/4)-f(a+(b-a)/2)+f(a+3(b-a)/4)\right\} .
\end{equation}

\bigskip 

\bigskip A third such approximation (figure 3) adds the area

\begin{eqnarray}
A_{3} &=&\frac{1}{8}(b-a)\{f(a+(b-a)/8)-f(a+(b-a)/4) \\
&&+f(a+3(b-a)/8)-f(a+(b-a)/2)+f(a+5(b-a)/8)  \notag \\
&&-f(a+3(b-a)/4)+f(a+7(b-a)/8)\}.  \notag
\end{eqnarray}

\bigskip 

\bigskip

Each new iteration n leads to 2$^{n}$-1 new triangles adding area A$_{n}$ to
further refine the area estimate. This procedure is continued indefinitely
to ``exhaust'' the remaining area, leading to the formula

\begin{equation}
\int_{a}^{b}f(x)dx=\sum_{n=1}^{\infty }A_{n}=(b-a)\sum_{n=1}^{\infty
}\sum_{m=1}^{2^{n}-1}(-1)^{m+1}2^{-n}f(a+m(b-a)/2^{n}).
\end{equation}

The special case for (5) when a = 0 is as follows:

\begin{equation}
\int_{0}^{b}f(x)dx=b\sum_{n=1}^{\infty
}\sum_{m=1}^{2^{n}-1}(-1)^{m+1}2^{-n}f(mb/2^{n}).
\end{equation}

The method of exhaustion will converge to the value of the integral at least
as fast as a geometric series, because when each new triangle is small
enough so that the local curvature between intersection points on f(x) is
slowly varying, it will exhaust more than 1/2 of the remaining area left
from the previous approximation.

There are an infinite number of possible variations on this procedure. For
example, each new iteration may involve a different number of triangles to
exhaust the remaining area, or each new triangle could intersect the
function at a location other than the midpoint.

\section{Algebraic Proof}

\bigskip An independent proof for (6) can be obtained starting with the
identity [2]

\begin{equation}
\frac{\sin a}{a}=\cos ^{2}(a/2)+\sum_{n=1}^{\infty }\sin
^{2}(a/2^{n+1})\prod_{m=1}^{n}\cos (a/2^{m})\ .
\end{equation}

\bigskip This expression can be rewritten as follows: 
\begin{equation}
\frac{\sin a}{a}=\sum_{n=1}^{\infty
}\sum_{m=1}^{2^{n}-1}(-1)^{m+1}2^{-n}\cos (ma/2^{n})\ ;
\end{equation}

\bigskip or

\begin{equation}
\frac{\sin ba}{a}=b\sum_{n=1}^{\infty
}\sum_{m=1}^{2^{n}-1}(-1)^{m+1}2^{-n}\cos (mba/2^{n})\ .
\end{equation}

\bigskip Since 
\begin{equation}
\frac{\sin ba}{a}=\frac{1}{2}\int_{-b}^{b}e^{-iax}dx,
\end{equation}

it follows that

\begin{equation}
\sum_{k=0}^{\infty }i^{k}\gamma _{k}\frac{d^{k}}{da^{k}}\left[ \sin (ba)/a%
\right] =\frac{1}{2}\int_{-b}^{b}\sum_{k=0}^{\infty }\gamma
_{k}x^{k}e^{-iax}dx.
\end{equation}

Defining 
\begin{equation}
G(x)=\sum_{k=0}^{\infty }\gamma _{k}x^{k},
\end{equation}

and noting that 
\begin{eqnarray}
\sum_{k=0}^{\infty }i^{k}\gamma _{k}\frac{d^{k}}{da^{k}}[\sin (ba)/a]
&=&b\sum_{n=1}^{\infty
}\sum_{m=1}^{2^{n}-1}(-1)^{m+1}2^{-n}\sum_{k=0}^{\infty }\gamma
_{2k}(mb/2^{n})^{2k}  \notag \\
&&\times \cos (mba/2^{n})  \notag \\
&&-ib\sum_{n=1}^{\infty
}\sum_{m=1}^{2^{n}-1}(-1)^{m+1}2^{-n}\sum_{k=0}^{\infty }\gamma
_{2k+1}(mb/2^{n})^{2k+1}  \notag \\
&&\times \sin (mba/2^{n}),
\end{eqnarray}

\bigskip

\bigskip it follows that

\begin{eqnarray}
\int_{-b}^{b}G(x)e^{-iax}dx &=&b\sum_{n=1}^{\infty
}\sum_{m=1}^{2^{n}-1}(-1)^{m+1}2^{-n}\left[ G(mb/2^{n})+G(-mb/2^{n})\right] 
\notag \\
&&\times \cos (mba/2^{n})  \notag \\
&&-ib\sum_{n=1}^{\infty }\sum_{m=1}^{2^{n}-1}(-1)^{m+1}2^{-n}\left[
G(mb/2^{n})-G(-mb/2^{n})\right]  \notag \\
&&\times \sin (mba/2^{n}).
\end{eqnarray}

Setting 
\begin{equation}
G(x)=H(x)e^{iax}
\end{equation}

leads to

\begin{equation}
\int_{-b}^{b}H(x)dx=b\sum_{n=1}^{\infty
}\sum_{m=1}^{2^{n}-1}(-1)^{m+1}2^{-n} \left[ H(mb/2^{n})+H(-mb/2^{n})\right]
.
\end{equation}

Note that H(x) can be in the form of a Fourier series which sums to the
following:

\bigskip

\begin{equation}
H(x)=\left\{ 
\begin{array}{c}
0;-b\leq x<0 \\ 
f(x);0\leq x<b
\end{array}
\right\} ,
\end{equation}

so that

\begin{equation}
\int_{0}^{b}f(x)dx=b\sum_{n=1}^{\infty
}\sum_{m=1}^{2^{n}-1}(-1)^{m+1}2^{-n}f(mb/2^{n}).
\end{equation}

\section{Applications}

The most immediate applications are series expressions for common functions
that might otherwise be very difficult to derive. Some examples are
presented below:

\bigskip 
\begin{equation}
\sin x=x\sum_{n=1}^{\infty }\sum_{m=1}^{2^{n}-1}(-1)^{m+1}2^{-n}\cos
(mx/2^{n})\ ;
\end{equation}

\begin{equation}
\cos x=1-x\sum_{n=1}^{\infty }\sum_{m=1}^{2^{n}-1}(-1)^{m+1}2^{-n}\sin
(mx/2^{n})\ ;
\end{equation}

\begin{equation}
\int_{0}^{b}\frac{\sin ax}{x}dx=\sum_{n=1}^{\infty }\sum_{m=1}^{2^{n}-1}%
\frac{(-1)^{m+1}}{m}\sin (mba/2^{n})\ ;
\end{equation}

\begin{equation}
e^{x}=1+x\sum_{n=1}^{\infty
}\sum_{m=1}^{2^{n}-1}(-1)^{m+1}2^{-n}e^{mx/2^{n}}\ ;
\end{equation}

\begin{equation}
\int_{0}^{b}e^{-ax^{2}}dx=b\sum_{n=1}^{\infty
}\sum_{m=1}^{2^{n}-1}(-1)^{m+1}2^{-n}e^{-a(mb)^{2}/4^{n}}\ ;
\end{equation}

\begin{eqnarray}
\ln x &=&\sum_{n=1}^{\infty }\sum_{m=1}^{2^{n}-1}\frac{(-1)^{m+1}(x-1)}{%
2^{n}+m(x-1)}; \\
(0 &<&x<\infty )  \notag
\end{eqnarray}

\begin{eqnarray}
p! &=&\int_{0}^{1}\left[ \ln (1/x)\right] ^{p}dx=\sum_{n=1}^{\infty
}\sum_{m=1}^{2^{n}-1}(-1)^{m+1}2^{-n}\left[ \ln (2^{n}/m)\right] ^{p}. \\
(0 &\leq &p<\infty )  \notag
\end{eqnarray}

In addition to the above expressions, this method can lead to new insights
into certain physical problems. One example involves the diffraction of
waves though a two-dimensional aperture in an infinite screen. This problem
has been solved exactly in integral form [3-5] so that the field at every
point depends on the Fourier transform F(k$_{x}$,k$_{y}$) of the aperture at
z=0:

\begin{equation}
\phi (x,y,z,t)=\frac{e^{-i\omega _{0}t}}{2\pi }\int_{-\infty }^{\infty
}\int_{-\infty }^{\infty }F(k_{x},k_{y})e^{-ik_{x}x}e^{-ik_{y}y}e^{iz\sqrt{%
k^{2}-k_{x}^{2}-k_{y}^{2}}}dk_{x}dk_{y}
\end{equation}

\bigskip Note that (26) satisfies the Helmholtz equation,

\begin{equation}
\nabla ^{2}\phi +k^{2}\phi =0
\end{equation}
everywhere (k = $\omega _{0}$/c) as well as the boundary conditions on the
screen.

Under suitable conditions, the extension of (5) to convergent improper
integrals can be made as follows, i.e.,

\bigskip 
\begin{equation}
\int_{0}^{\infty }f(x)dx=b\sum_{n=1}^{\infty
}\sum_{m=1}^{2^{n}-1}\sum_{p=0}^{\infty }(-1)^{m+1}2^{-n}f(pb+mb/2^{n}).
\end{equation}

Equation (28) is a result of breaking the integral into a series of definite
integrals on [pb,(p+1)b]. The result is valid when each integral 
\begin{equation*}
\int_{pb}^{(p+1)b}f(x)dx
\end{equation*}
exists and when

\bigskip 
\begin{equation}
\int_{0}^{B}f(x)dx
\end{equation}
tends to a finite limit L as B$\rightarrow \infty .$

With regard to (26), however, it can be noted that propagating waves will
only occur when the transverse wavenumber (k$_{x}^{2}$+k$_{y}^{2}$)$^{1/2}$
is lower than the cutoff wavenumber, so in the far field the integral can be
evaluated to $\pm $k with good accuracy as follows (assuming that F(k$_{x}$,k%
$_{y}$) is an even function with respect to both k$_{x}$ and k$_{y}$):

\bigskip

\begin{eqnarray}
\phi (x,y,z,t) &\simeq &\frac{2k^{2}}{\pi }e^{-i\omega
_{0}t}\sum_{n=1}^{\infty }\sum_{m=1}^{2^{n}-1}\sum_{p=1}^{\infty
}\sum_{q=1}^{2^{p}-1}(-1)^{m+q}2^{-n-p}F\left( mk/2^{n},qk/2^{p}\right) 
\notag \\
&&\times \cos (mxk/2^{n})\cos (qyk/2^{p})e^{izk\sqrt{%
1-m^{2}/4^{n}-q^{2}/4^{p}}}
\end{eqnarray}

It can be clearly seen from (30) that the resulting field is due to a
summation of an infinite number of plane waves, each propagating at a
different phase velocity (and hence a different group velocity) based on its
value of k$_{x}$ and k$_{y}$. Thus, the process of diffraction leads to a
continuous distribution of group velocities, having an amplitude
distribution governed by the Fourier transform of the aperture. This insight
becomes particularly clear upon application of (5) to expand (26). Such a
result has an effect on Doppler shifts both in acoustics and in
electromagnetic wave propagation, effects that are not otherwise apparent
without the use of the generalized method of exhaustion.

\bigskip

\bigskip

\section{References}

\begin{enumerate}
\item  Simmons, G. F., ``Calculus with Analytic Geometry''. New York:
McGraw-Hill, 1996, p. 190.

\item  Ruffa, A.A., \textit{A Series for sinx/x}. Math. Mag. (to appear
October 2000)

\item  Brekhovskikh, L.M., ``Waves in Layered Media''. New York: Academic
Press, 1960, p. 100.

\item  Born, M. and Wolf, E., ``Principles of Optics (Seventh Edition)''.
Cambridge University Press, 1999, p. 640.

\item  Gaunaurd, G.C. and Uberall, H., \textit{Acoustics of Finite Beams}.
J. Acoust. Soc. Am. 63, 5, (1978).
\end{enumerate}

\end{document}